\documentclass[11pt,a4paper]{article}

\usepackage[T1]{fontenc}
\usepackage{lmodern}
\usepackage[margin=27mm]{geometry}

\usepackage{amsmath,amssymb,amsthm,mathtools}
\usepackage{booktabs}
\usepackage{array}
\usepackage{enumitem}
\usepackage{microtype}
\usepackage{tikz}
\usetikzlibrary{arrows.meta,positioning,fit,calc}

\usepackage{xurl}
\usepackage[hidelinks]{hyperref}
\usepackage[nameinlink,noabbrev,capitalise]{cleveref}

\usepackage[
    backend=biber,
    style=numeric,
    sorting=none
]{biblatex}

\theoremstyle{definition}
\newtheorem{definition}{Definition}
\newtheorem{principle}{Principle}
\newtheorem{invariant}{Invariant}
\theoremstyle{remark}
\newtheorem{remark}{Remark}

\crefname{principle}{Principle}{Principles}
\crefname{invariant}{Invariant}{Invariants}
\crefname{definition}{Definition}{Definitions}
\crefname{remark}{Remark}{Remarks}

\newcommand{\State}{\mathcal{S}}
\newcommand{\Input}{\mathcal{X}}
\newcommand{\Output}{\mathcal{O}}
\newcommand{\Random}{\mathcal{R}}
\newcommand{\Term}{\State_{\top}}
\newcommand{\Traj}{\mathcal{T}}
\newcommand{\InputContract}{\mathsf{IC}}
\newcommand{\Exp}{\mathbb{E}}
\newcommand{\Prob}{\mathbb{P}}
\newcommand{\code}[1]{\texttt{#1}}

\tikzset{
    arch/.style={
        draw,
        rounded corners=2pt,
        align=center,
        minimum width=29mm,
        minimum height=10mm,
        font=\small
    },
    small/.style={
        draw,
        rounded corners=2pt,
        align=center,
        minimum width=18mm,
        minimum height=8mm,
        font=\small
    },
    lifecycle/.style={
        draw,
        rounded corners=4pt,
        align=center,
        minimum width=22mm,
        minimum height=8mm,
        font=\small\scshape
    },
    arrow/.style={
        -{Latex[length=2.2mm]},
        semithick
    },
    dashedarrow/.style={
        -{Latex[length=2.2mm]},
        semithick,
        dashed
    },
    lbl/.style={font=\scriptsize}
}

\title{
    Separating Statistical Inference from Experimental Runtime:\\
    An Architecture for Sequential Scientific Analysis
}

\author{
    Alexander Stepanov\\
    \small Sobolev Institute of Mathematics
    \and
    Andrej Novikov\\
    \small Novikov Laboratories LLC
}

\date{}

\begin{document}

\maketitle

\begin{abstract}

Sequential statistical procedures are increasingly embedded in scientific
systems in which observations arrive incrementally from heterogeneous
sources. In such systems the statistical update is only one part of the
computation: observations must also be acquired, aligned across inputs,
ordered, delivered, persisted, and replayed, and results must be passed to
systems that display them or act on them. Sequential methods are
nevertheless missing from the general software toolkit. They run inside
proprietary experimentation platforms and commercial statistical
packages, or exist as offline functions over complete data, and the batch
\code{fit} interface that dominates machine-learning libraries has no
place for a procedure that listens, keeps state, and decides when to
stop.

We propose an architecture that separates statistical inference from data
topology, execution runtime, run lifecycle, persistence, and external
action. Its central abstraction is the \emph{Sequential Procedure}: a
state-transition function that maps an inferential state, a valid
statistical input, and a statistical random state to a new state and a
scientific output. The procedure does not define how observations are
acquired, paired, scheduled, stored, or acted upon. We specify the
architecture semantically rather than through an object model, derive
seven conformance invariants that any implementation can be tested
against (among them synchronous/asynchronous, checkpoint, replay, and
batch equivalence), and describe a dependency-free reference
implementation in which all seven hold bitwise. Lorden's 2-SPRT for the
Kiefer--Weiss problem serves as the reference procedure; the same
contract accommodates sequential estimators, Bayesian filters,
change-point detectors, and stateful predictive models.

\end{abstract}

\noindent
\textbf{Keywords:}
sequential inference,
sequential statistics,
scientific computing,
streaming data,
experimental runtime,
checkpointing,
provenance,
data topology,
Kiefer--Weiss,
stateful inference

\section{Introduction}
\label{sec:introduction}

Many scientific analyses are inherently sequential. Observations are
collected over time, the available information changes with every
observation, and the current inferential state may determine what is
measured next or whether measurement continues at all. Sequential
hypothesis testing~\cite{wald1945,wald1947}, change-point
detection~\cite{page1954,poor2009}, Bayesian filtering and state-space
inference~\cite{kalman1960,sarkka2013}, and sequential prediction are all
of this kind.

The mathematical definition of such a procedure is usually small: given
a state and a new observation, produce a new state and an output. A
working scientific implementation has to solve a much larger problem.
Observations originate from heterogeneous sources, arrive asynchronously,
need to be aligned with one another, may be delivered in batches,
duplicated, or reordered, and have to be persisted so that an experiment
that runs for days can be interrupted and resumed.

\paragraph{Fit is not listen.}
The interface that dominates machine-learning software does not have this
shape. A \code{fit} method receives a complete data set and returns a
fitted object~\cite{buitinck2013}. A sequential procedure has to
\emph{listen}: it receives one observation, updates a state, and may
announce that it has seen enough. Incremental variants such as
\code{partial\_fit} process data in chunks, but they treat samples as
exchangeable. They carry no step index, no stopping rule, and no state
that means anything as a point in a sequence. Even models that are
sequential by construction are commonly used through the batch interface.
A widely used hidden Markov model library offers estimation, decoding,
and scoring of whole sequences~\cite{hmmlearn}; the forward recursion
that updates a belief after each observation is not part of its public
interface. Recurrent networks are often trained and applied on fixed
windows, with the hidden state reset between windows. In both cases the
state that makes the model sequential is discarded at the boundary of
every call.

\paragraph{Where sequential methods live today.}
Sequential methods with decision semantics are in production use, yet
they are outside the general toolkit (\cref{tab:landscape}). Commercial
experimentation platforms run mixture SPRTs and anytime-valid confidence
sequences on live data~\cite{johari2017,johari2022,zhao2019,lindon2022}.
Commercial statistical software and a validated R package provide
group-sequential designs for clinical
trials~\cite{sas,rpact,wassmer2025}. The former are proprietary and
embedded in a platform. The latter operate on assembled data sets at
interim analyses. Open implementations of fully sequential tests are
mostly R functions that take the complete vector of observations and
return a decision~\cite{schnuerch2020,pramanik2021,grunwald2024}. That
is of limited use in application engineering. R is an environment for
analysis and is rarely the language in which a deployed system is
written: in the 2025 Stack Overflow survey, 3.3\% of professional
developers reported extensive development work in R, against 54.8\% for
Python~\cite{stackoverflow2025}. Using an R function from a service
means operating a second runtime behind a bridge, or re-implementing the
method, and that cost is paid again for every method adopted. Open
incremental libraries exist for change detection and drift
alarms~\cite{romano2026,montiel2021}; their state is mutable and
internal, and the stopping rule is left to the caller. We found no
packaged implementation of Lorden's 2-SPRT. The available code for the
Kiefer--Weiss problem, including our own, computes optimal designs
offline~\cite{novikov2022a,novikov2022b,novikov2023}.

The gap is therefore not in theory and not in demand. What is missing is
a generic, open component model for sequential procedures, as opposed to
one more implementation of one method inside one system. Numerical
computing was once in the same position. The algorithms existed as
Fortran routines, and what general-purpose libraries such as
SciPy~\cite{virtanen2020} added was a uniform and accessible interface
to them, with no new numerical methods. That step is what turned
specialist routines into a general toolkit, and sequential inference has
not yet taken it.

When such methods are implemented ad hoc, the surrounding concerns end up
inside the statistical method.
The implementation of a test or a filter then depends on a particular
event loop, message broker, database, sensor API, file format, or
plotting library. None of these belongs to the mathematical definition of
the procedure, yet they become inseparable from its software
representation. The consequences are practical. The method cannot be
reused in another laboratory without its infrastructure, it cannot be
validated in isolation, and a run cannot be reproduced without
reproducing the environment that happened to execute it.

This paper proposes a separation, stated as a single principle:

\begin{quote}
A sequential statistical procedure should define how valid observations
update inferential state and produce scientific outputs. It should not
define how observations are acquired, scheduled, transported, persisted,
or acted upon.
\end{quote}

Two aspects of this separation deserve emphasis because they are where
existing designs most often leak.

The first is \emph{topology}. A sequential scientific system is not a
single stream-processing function. A procedure may receive two
independent samples, paired observations, or several heterogeneous
inputs with different acquisition mechanisms and different clocks.
Whether $A_3$ is analyzed together with $B_3$, with the $B$ nearest in
time, or on its own is a property of the experimental design. It is not
a property of the likelihood update, and it should not be reconstructed
inside it.

The second is \emph{action}. A procedure produces scientific outputs;
other components may use those outputs for visualization, logging,
experimental design, alerting, or control. We summarize this as

\begin{quote}
\textbf{Procedure infers; consumers act.}
\end{quote}

\paragraph{Contributions.}
The paper makes four contributions.
\begin{enumerate}[label=(\roman*),nosep]
    \item A semantic definition of a Sequential Procedure as a
    state-transition function with an explicit statistical random state,
    partial multi-input structure, and absorbing terminal states
    (\cref{sec:procedure}).
    \item A layered architecture around that definition, with precise
    responsibilities for adapters, input contracts, data topology,
    runtime, run, and consumers (\cref{sec:architecture,sec:inputs,%
    sec:runtime,sec:run,sec:consumers}).
    \item Seven conformance invariants that follow from the
    architectural principles and can be checked mechanically for any
    implementation (\cref{sec:invariants}).
    \item A reference implementation with no dependencies in which the
    invariants are executable tests, and a worked reference procedure
    for the Kiefer--Weiss problem
    (\cref{sec:kw,sec:implementation}).
\end{enumerate}

\section{Background and Related Work}
\label{sec:related}

\paragraph{Sequential analysis.}
The statistical side of the problem is classical. Wald's sequential
probability ratio test (SPRT)~\cite{wald1945,wald1947} and its
optimality~\cite{waldwolfowitz1948}, Page's cumulative-sum
scheme~\cite{page1954}, and the theory built on
them~\cite{siegmund1985,poor2009,tartakovsky2014} define procedures by a
recursion on a low-dimensional statistic and a stopping rule. The more
recent literature on anytime-valid inference~\cite{ramdas2023,howard2021}
does the same with e-processes and confidence sequences. Filtering and
state-space
methods~\cite{kalman1960,rabiner1989,doucet2001,sarkka2013,bishop2006,%
murphy2022} and recurrent predictive models~\cite{hochreiter1997} share
that recursive shape. In all of these the published object is a
transition rule. What the literature does not prescribe, reasonably, is
how the rule should be embedded in a running system.

\paragraph{Software for sequential inference.}
\Cref{tab:landscape} summarizes where sequential inference is available
as software. The table was compiled in October 2026 from package
documentation and published descriptions. It is not exhaustive, and for
proprietary systems it reflects only what their authors have published.

\begin{table}[t]
\centering
\small
\caption{
Where sequential inference is available as software. ``Consumes'' is the
unit of data handed to the method in one call. ``State'' is what the
caller can keep between calls in order to continue.
}
\label{tab:landscape}
\begin{tabular}{
    >{\raggedright\arraybackslash}p{0.27\linewidth}
    >{\raggedright\arraybackslash}p{0.14\linewidth}
    >{\raggedright\arraybackslash}p{0.15\linewidth}
    >{\raggedright\arraybackslash}p{0.17\linewidth}
    >{\raggedright\arraybackslash}p{0.12\linewidth}}
\toprule
\textbf{System} & \textbf{Access} & \textbf{Consumes} & \textbf{State} & \textbf{Stopping rule} \\
\midrule
Optimizely, Statsig \cite{johari2017,zhao2019}
    & proprietary platform & live metrics & not public & yes \\
Netflix \cite{lindon2022}
    & internal & live counts & not public & yes \\
SAS SEQDESIGN, SEQTEST \cite{sas}
    & commercial & interim stage & boundary data set & yes \\
rpact \cite{rpact}
    & open, R & interim stage & R objects & yes \\
sprtt, MSPRT \cite{schnuerch2020,pramanik2021}
    & open, R & full vector & none & yes \\
safestats \cite{grunwald2024}
    & open, R & full vector & none & yes \\
\addlinespace
confseq \cite{howard2021}
    & open, Python & running statistics & held by caller & left to caller \\
Confidence \cite{spotifyconfidence}
    & open, Python & data frame per look & none & bounds only \\
focus, focus-cpt \cite{romano2026}
    & open, R and Python & one observation & mutable, internal & left to caller \\
River drift detectors \cite{montiel2021}
    & open, Python & one observation & mutable, internal & alarm flag \\
Spark \code{StreamingTest} \cite{spark}
    & open, JVM & micro-batch & not exposed & none \\
\addlinespace
\code{partial\_fit} \cite{buitinck2013}
    & open, Python & mini-batch & model parameters & none \\
hmmlearn \cite{hmmlearn}
    & open, Python & whole sequence & none & --- \\
Kiefer--Weiss code \cite{novikov2022a}
    & open, R scripts & design parameters & --- & plan only \\
\addlinespace
This work
    & open, Python & one statistical input & explicit, serializable, with random state & terminal state \\
\bottomrule
\end{tabular}
\end{table}

Three groups can be distinguished. The first has decision semantics and
is production-grade, but it is proprietary, commercially licensed, or
tied to interim analyses of assembled data. The mixture SPRT of Johari
et~al.~\cite{johari2017,johari2022} and the staged-rollout tests of Zhao
et~al.~\cite{zhao2019} are examples of the former; group-sequential
software for clinical trials~\cite{sas,rpact,wassmer2025} is an example
of the latter. SAS's analysis procedure is notable because it passes a
boundary data set from one interim stage to the next, which is a
documented form of resumable state at the granularity of stages. The
second group is open and has decision semantics, but consumes complete
data vectors and keeps nothing between
calls~\cite{schnuerch2020,pramanik2021,grunwald2024}. The third group is
open and incremental, but leaves the state implicit and the decision to
the caller: boundary functions for confidence
sequences~\cite{howard2021}, online change-point
detectors~\cite{romano2026}, drift detectors~\cite{montiel2021}, and a
streaming two-sample test that reports fixed-horizon $p$-values without
a stopping rule~\cite{spark}. None of the systems we examined offers a
sequential test as a component that consumes one observation at a time,
exposes a terminal decision, and keeps its inferential and random state
in a form that can be persisted and resumed. For the Kiefer--Weiss
problem, earlier work co-authored by one of
us~\cite{novikov2022a,novikov2023} and its
companion~\cite{novikov2022b} compute optimal tests by backward
induction. The accompanying code produces the plan, a table of
continuation intervals, and does not execute it on a stream
(\cref{sec:kw}).

\paragraph{Incremental learning interfaces.}
Libraries for online machine learning expose incremental update methods
on mutable model objects~\cite{bifet2010,buitinck2013,montiel2021,%
day2020}. This is close in spirit to a Sequential Procedure. The
differences lie in what is made explicit: these interfaces leave
ordering, persistence of randomness, multi-input topology, stopping, and
the provenance of a run to the caller.

\paragraph{Stream processing.}
Distributed stream processors solved many of the execution problems we
discuss: event time versus processing time, windowing, out-of-order
data, and consistent snapshots of operator
state~\cite{akidau2015,carbone2015,carbone2017,kreps2011,chandy1985,%
elnozahy2002}. Their determinacy results go back to Kahn's process
networks~\cite{kahn1974}. Summingbird~\cite{boykin2014} showed that pure
aggregators can be executed unchanged by a batch and a streaming runtime,
which is a direct precedent for batch transparency. It relies on
associativity, which a procedure with a stopping rule does not have.
These systems treat operator state as opaque and optimize for throughput
and fault tolerance at scale. They do not say what a statistical random
state is, when two state trajectories count as equal, or what a terminal
decision means for acquisition. Our architecture borrows their
vocabulary and is meant to be hosted by such systems as readily as by a
for-loop; it specifies the part they leave open.

\paragraph{Experiment orchestration.}
Frameworks for instrument control separate the description of an
experiment from the engine that executes it and from the consumers of
its data~\cite{allan2019}. This is the closest precedent for our
separation of acquisition and consumers. Such frameworks do not define a
contract for the statistical procedure that sits between the two.

\paragraph{Reproducibility and provenance.}
Guidelines for reproducible computational
research~\cite{sandve2013,wilkinson2016} and provenance
models~\cite{moreau2013} ask that the lineage of a result be recorded.
For sequential experiments the lineage includes things that batch
analyses do not have: the order in which observations were applied, what
was dropped under load, and the random state at every interruption.

\paragraph{Software architecture and novelty.}
The separation itself is an instance of information
hiding~\cite{parnas1972} and of the ports-and-adapters
style~\cite{cockburn2005}. Formally, a Sequential Procedure is a Mealy
machine~\cite{mealy1955} with a possibly infinite state space. Neither
the automaton nor the idea of separating pure logic from runtime and
persistence is new. What we add is the application of that pattern to
sequential procedures with decision semantics, and the parts it then
requires: a statistical random state that is carried in the contract, a
topology layer that relates multiple inputs before inference, and
conformance invariants against which any implementation can be tested.
As with general-purpose scientific libraries~\cite{virtanen2020}, the
value of a common contract does not depend on the methods behind it
being new. Existing kernels, including those written in R or C++, can be
placed behind it without being rewritten.

\section{Architectural Principles}
\label{sec:principles}

The architecture rests on seven principles. Each is used later to derive
a testable invariant (\cref{sec:invariants}).

\begin{principle}[Inference--execution separation]
\label{pr:execution}
The state transition of a procedure is independent of the mechanism that
delivers its inputs.
\end{principle}

\begin{principle}[Topology--inference separation]
\label{pr:topology}
Joining, pairing, alignment, filtering, mapping, and windowing of data
sources occur outside the procedure.
\end{principle}

\begin{principle}[Logical inputs]
\label{pr:logical}
A procedure consumes logical statistical inputs, not physical sources.
Physical acquisition is connected through adapters.
\end{principle}

\begin{principle}[Single-writer inference]
\label{pr:single-writer}
For one procedure instance, state transitions are serialized according
to an explicit ordering policy.
\end{principle}

\begin{principle}[Persistent state]
\label{pr:persistent}
The state required to continue a run is distinct from, and independent
of, the disposable execution state of the runtime.
\end{principle}

\begin{principle}[Consumer non-interference]
\label{pr:consumers}
Consumers observe outputs. They cannot modify inferential state, and
their failures are not inference failures.
\end{principle}

\begin{principle}[Semantic portability]
\label{pr:portability}
The architecture is defined by semantic contracts, not by a particular
programming-language object model.
\end{principle}

\section{System Overview}
\label{sec:architecture}

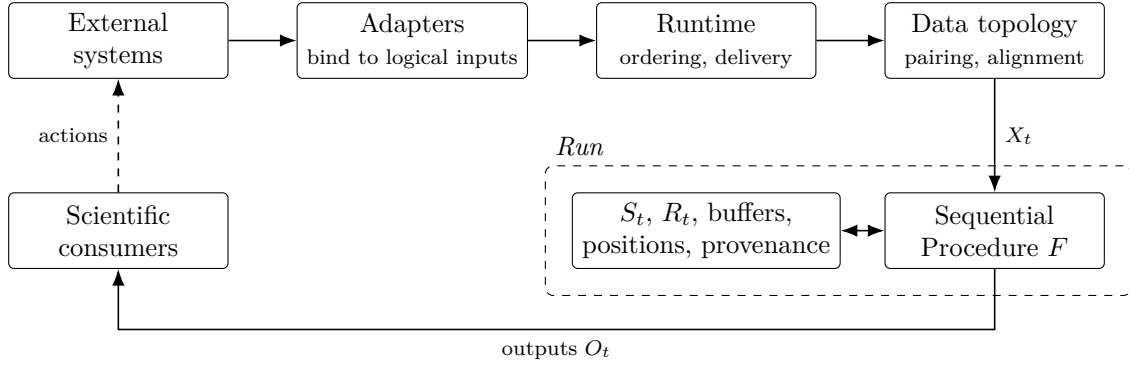
\begin{figure}[t]
\centering
\begin{tikzpicture}[node distance=15mm and 9mm]

\node[arch] (world) {External\\systems};
\node[arch, right=of world] (adapter)
    {Adapters\\{\scriptsize bind to logical inputs}};
\node[arch, right=of adapter] (runtime)
    {Runtime\\{\scriptsize ordering, delivery}};
\node[arch, right=of runtime] (topology)
    {Data topology\\{\scriptsize pairing, alignment}};

\node[arch, below=of topology] (procedure) {Sequential\\Procedure $F$};
\node[arch, below=of runtime] (store)
    {$S_t$, $R_t$, buffers,\\positions, provenance};
\node[arch, below=of world] (consumer) {Scientific\\consumers};

\node[draw, dashed, rounded corners=3pt, inner sep=3.5mm,
      fit=(procedure)(store)] (run) {};
\node[lbl, anchor=south west, font=\small\itshape]
    at (run.north west) {Run};

\draw[arrow] (world) -- (adapter);
\draw[arrow] (adapter) -- (runtime);
\draw[arrow] (runtime) -- (topology);
\draw[arrow] (topology) -- node[right,lbl] {$X_t$} (procedure);
\draw[arrow, {Latex[length=2.2mm]}-{Latex[length=2.2mm]}]
    (procedure) -- (store);
\draw[arrow] (procedure.south) -- ++(0,-8mm)
    -| node[pos=.25,below,lbl] {outputs $O_t$} (consumer.south);
\draw[dashedarrow] (consumer) -- node[left,lbl] {actions} (world);

\end{tikzpicture}
\caption{
The data path. The run is not a stage on the path: it is the context
that owns everything which must outlive the executing process. Actions
return to the external world, never directly to the procedure.
}
\label{fig:architecture}
\end{figure}

\Cref{fig:architecture} shows the data path. Observations leave an
external system, are bound to logical inputs by an adapter, are ordered
and deduplicated per source by the runtime, and are related to one
another by the data topology. What crosses into the run is a
\emph{statistical input} $X_t$. The procedure maps it to an output $O_t$,
which is handed to consumers.

The layers have the following responsibilities.

\begin{description}[leftmargin=2em,style=nextline,itemsep=2pt]
    \item[Adapters] connect external systems to logical inputs. They
    convert formats and units and attach source identity, sequence
    numbers, keys, and event times.
    \item[Logical inputs] carry a contract that states what a valid
    observation is (\cref{sec:contracts}).
    \item[Runtime] decides how observations are acquired and scheduled,
    and in which per-source order they are delivered
    (\cref{sec:runtime}). It holds no state that needs to survive.
    \item[Data topology] relates inputs to one another and emits
    statistical inputs (\cref{sec:topology}).
    \item[Run] is one concrete execution. It owns the inferential state,
    the statistical random state, topology buffers, source positions,
    lifecycle, and provenance, and it is the single writer of the
    inferential state (\cref{sec:run}).
    \item[Sequential Procedure] performs the state transition
    (\cref{sec:procedure}).
    \item[Scientific consumers] use outputs for observation, reporting,
    design, or action (\cref{sec:consumers}).
\end{description}

Per-source ordering precedes topology because positional pairing is only
meaningful on ordered inputs. Serialization of state transitions follows
topology because that is where statistical inputs come into existence.

Not every application instantiates every layer. A minimal offline
analysis reduces to a procedure folded over an iterator. A laboratory
system instantiates all of them. The procedure is the same object in
both.

\section{Sequential Procedure Semantics}
\label{sec:procedure}

\begin{definition}[Sequential Procedure]
\label{def:procedure}
A sequential procedure is a tuple
$(\State,\Input,\Output,\Random,s_0,F,\Term)$ consisting of a state space
$\State$, an input space $\Input$, an output space $\Output$, a space
$\Random$ of statistical random states, an initial state
$s_0\in\State$, a transition function
\[
F:\State\times\Input\times\Random\to\State\times\Output\times\Random ,
\]
and a set $\Term\subseteq\State$ of terminal states.
\end{definition}

Given an initial state $S_0=s_0$, an initial random state $R_0$, and
statistical inputs $X_1,X_2,\ldots$, the procedure evolves as
\begin{equation}
(S_t,\,O_t,\,R_t) \;=\; F(S_{t-1},\,X_t,\,R_{t-1}),
\qquad t=1,2,\ldots
\label{eq:transition}
\end{equation}
so that $S_t$ and $O_t$ are functions of $X_{1:t}$. For a deterministic
procedure $\Random$ is a single point and \eqref{eq:transition} reduces
to $(S_t,O_t)=f(S_{t-1},X_t)$.

\begin{definition}[Trajectory]
\label{def:trajectory}
For an input history $x_{1:T}$, the trajectory is
\[
\Traj(x_{1:T};s_0,r_0) \;=\; \bigl((S_t,O_t,R_t)\bigr)_{t=1}^{\tau},
\qquad
\tau=\min\bigl(T,\;\inf\{t: S_t\in\Term\}\bigr).
\]
\end{definition}

The trajectory is the reference semantics of a run. Every claim of the
form ``mechanism $A$ is equivalent to mechanism $B$'' in this paper
means that both produce the same trajectory, under a stated notion of
equality (\cref{sec:numerical}).

The definition is deliberately broader than hypothesis testing. The
state may be a test statistic, an estimator, a posterior distribution, a
latent-state belief, or the hidden state of a predictive model;
\cref{tab:states} gives examples. The architecture imposes no universal
representation on $\State$. It requires only that a state can be encoded
to a value and decoded again without loss.

\begin{table}[t]
\centering
\caption{Examples of inferential state.}
\label{tab:states}
\begin{tabular}{lll}
\toprule
\textbf{Procedure} & \textbf{State $S_t$} & \textbf{Terminal} \\
\midrule
Likelihood-ratio test & $(t,\ \log\Lambda_t,\ d_t)$ & decision $d_t$ taken \\
CUSUM detector & $(t,\ W_t)$ & alarm raised \\
Kalman filter & $(\hat z_t,\ P_t)$ & never \\
Bayesian filter & $p(z_t\mid x_{1:t})$ or a particle cloud & never \\
LSTM predictor & $(h_t,\ c_t)$ & never \\
\bottomrule
\end{tabular}
\end{table}

\subsection{Statistical inputs}
\label{sec:multiple-inputs}

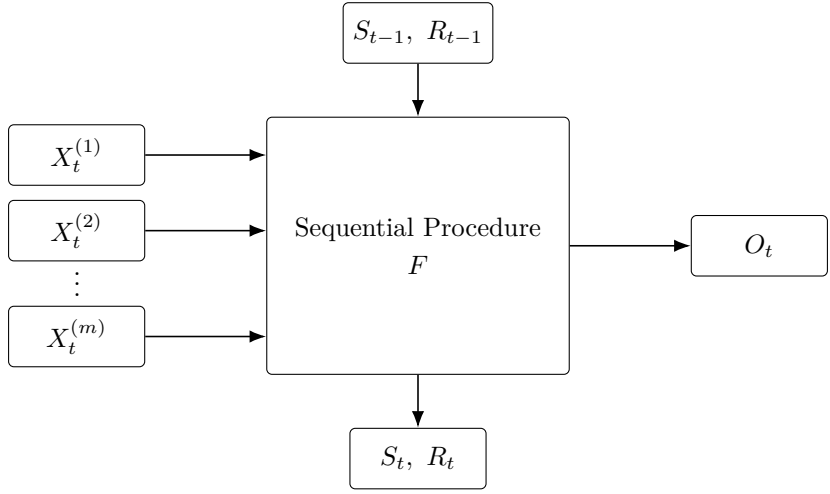
\begin{figure}[t]
\centering
\begin{tikzpicture}[node distance=4mm and 16mm]

\node[arch, minimum width=40mm, minimum height=34mm] (proc)
    {Sequential Procedure\\[2pt]$F$};

\node[small, anchor=east] (x1) at ($(proc.west)+(-16mm,12mm)$) {$X_t^{(1)}$};
\node[small, anchor=east] (x2) at ($(proc.west)+(-16mm,2mm)$) {$X_t^{(2)}$};
\node[font=\small] at ($(proc.west)+(-25mm,-4mm)$) {$\vdots$};
\node[small, anchor=east] (xm) at ($(proc.west)+(-16mm,-12mm)$) {$X_t^{(m)}$};

\node[small, above=7mm of proc] (prev) {$S_{t-1},\ R_{t-1}$};
\node[small, below=7mm of proc] (next) {$S_t,\ R_t$};
\node[small, right=of proc] (out) {$O_t$};

\draw[arrow] (x1.east) -- (x1.east -| proc.west);
\draw[arrow] (x2.east) -- (x2.east -| proc.west);
\draw[arrow] (xm.east) -- (xm.east -| proc.west);
\draw[arrow] (prev) -- (proc);
\draw[arrow] (proc) -- (next);
\draw[arrow] (proc) -- (out);

\end{tikzpicture}
\caption{
A procedure consumes logically distinct statistical inputs, any
non-empty subset of which may be present at a step. It receives
statistical inputs, not physical streams.
}
\label{fig:procedure}
\end{figure}

Let a procedure declare $m$ logical inputs with value spaces
$\Input^{(1)},\ldots,\Input^{(m)}$. Requiring every step to carry all of
them would exclude the most common two-sample situation, in which
observations of independent samples arrive one at a time and in no fixed
relation. We therefore define a statistical input as a partial
assignment
\begin{equation}
X_t \;=\; \bigl(X_t^{(j)}\bigr)_{j\in J_t},
\qquad \emptyset\neq J_t\subseteq\{1,\ldots,m\}.
\label{eq:multi-input}
\end{equation}
A procedure is \emph{joint} if it requires $J_t=\{1,\ldots,m\}$ at every
step, as a paired test does. Otherwise it accepts any non-empty subset.
The inputs stay logically distinct in either case
(\cref{fig:procedure}).

This avoids reducing the system to one undifferentiated event stream.
A single stream would leave the relation between inputs implicit and
force the procedure to reconstruct it, which is exactly the coupling
\cref{pr:topology} rules out.

\subsection{Sample structure and observation dimensionality}

The number of samples and the dimensionality of an observation are
independent axes. A one-sample multivariate procedure has
$X_t\in\mathbb{R}^d$; a two-sample scalar procedure has
$X_t^{(1)},X_t^{(2)}\in\mathbb{R}$; a two-sample multivariate procedure
has $X_t^{(1)},X_t^{(2)}\in\mathbb{R}^d$. Sample structure, observation
structure, and alignment are therefore expressed as contracts and
metadata. The architecture needs no separate type hierarchies for one-,
two-, and $k$-sample procedures.

\subsection{Terminal states}
\label{sec:terminal}

Terminal states are absorbing: once $S_t\in\Term$, no further transition
is applied and later statistical inputs do not alter the trajectory.
Making this part of the definition matters for two reasons. It fixes the
meaning of ``the same history'' when a decision falls in the middle of a
batch (\cref{sec:batching}). And it separates the statistical event
(a decision has been taken) from the operational one (acquisition
stops), which the runtime may schedule later (\cref{sec:lifecycle}).

\subsection{Randomized procedures}
\label{sec:random}

Three kinds of randomness occur in a sequential experiment and must not
share a generator.

\begin{enumerate}[nosep]
    \item \emph{Statistical randomness} is drawn by $F$ itself, as in a
    particle filter or a randomized test. Its state $R_t$ is part of the
    trajectory.
    \item \emph{Simulation randomness} generates synthetic observations.
    It belongs to a source.
    \item \emph{Execution randomness} arises from scheduling, retries,
    or jitter. It belongs to the runtime and must not influence the
    trajectory.
\end{enumerate}

Because $R_t$ appears in \eqref{eq:transition}, reproducing a run that
was interrupted at step $k$ requires $R_k$. The initial seed is not
enough unless the whole history is replayed. A checkpoint therefore
stores the current random state, and the transition commits $S_t$ and
$R_t$ together. For a fixed implementation and numerical contract,
$(x_{1:T},s_0,r_0)$ determines $\Traj(x_{1:T};s_0,r_0)$.

\section{Inputs: Contracts and Topology}
\label{sec:inputs}

\subsection{Input contracts}
\label{sec:contracts}

Every logical input $j$ carries a contract $\InputContract_j$. A contract
specifies the name and semantic identity of the input, its data type and
shape, its domain, its unit, its missing-value policy, and any metadata
required for ordering, alignment, or provenance. A value is valid if
\begin{equation}
x \models \InputContract_j .
\end{equation}

Validation happens before the transition. An input that fails its
contract never reaches $F$, so a contract violation cannot leave a
partially updated state. Because the contract belongs to the logical
input, a physical source can be replaced by any other source that
satisfies the same contract without touching the procedure
(\cref{pr:logical}).

\subsection{Data topology}
\label{sec:topology}

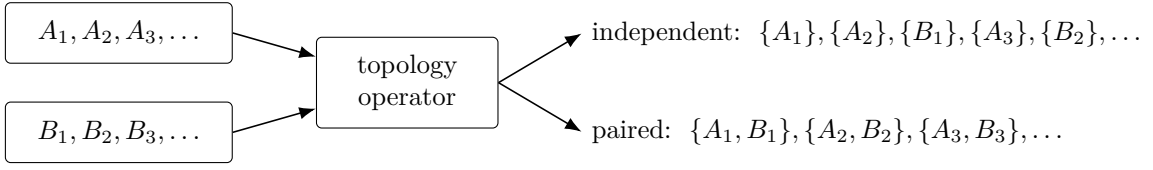
\begin{figure}[t]
\centering
\begin{tikzpicture}[node distance=5mm and 12mm]

\node[small, minimum width=30mm] (a) {$A_1,A_2,A_3,\ldots$};
\node[small, minimum width=30mm, below=of a] (b) {$B_1,B_2,B_3,\ldots$};
\coordinate (mid) at ($(a.east)!0.5!(b.east)$);
\node[small, minimum width=24mm, minimum height=12mm, anchor=west]
    (op) at ($(mid)+(11mm,0)$) {topology\\operator};

\node[font=\small, anchor=west, align=left] (ind) at ($(op.east)+(11mm,6.5mm)$)
    {independent:\ \ $\{A_1\},\{A_2\},\{B_1\},\{A_3\},\{B_2\},\ldots$};
\node[font=\small, anchor=west, align=left] (pair) at ($(op.east)+(11mm,-6.5mm)$)
    {paired:\ \ $\{A_1,B_1\},\{A_2,B_2\},\{A_3,B_3\},\ldots$};

\draw[arrow] (a.east) -- (op);
\draw[arrow] (b.east) -- (op);
\draw[arrow] (op.east) -- (ind.west);
\draw[arrow] (op.east) -- (pair.west);

\end{tikzpicture}
\caption{
The same two sources under two topologies. Pairing and alignment happen
before statistical inference, not inside it.
}
\label{fig:topology}
\end{figure}

Data topology specifies how logical inputs are related before they cross
the statistical boundary. Consider two sources $A$ and $B$ that produce
$A_1,A_2,\ldots$ and $B_1,B_2,\ldots$ respectively. For independent
samples each observation is a statistical input of its own. For a paired
analysis a topology operator constructs the pairs
(\cref{fig:topology}). The two analyses use the same sources and
different topologies. A paired test, moreover, needs no procedure of its
own: it is a one-sample test preceded by pairing and a difference map.

Topology operations include positional pairing, key-based joins,
timestamp alignment with a tolerance, windowing, filtering, and
stateless mapping. Two properties of these operators have consequences
for the rest of the architecture.

\paragraph{Topology operators have state.}
A pairing operator that has received $A_4$ but not yet $B_4$ holds an
observation that has been consumed from its source and has not reached
the procedure. If the run is interrupted at that moment, $A_4$ exists
only in the operator's buffer. Topology state is therefore part of what
must be persisted (\cref{sec:persistence}), and observations that an
operator discards, for example because no partner arrived within the
tolerance, change the observation history and must be counted.

\paragraph{Some operators are confluent.}
\begin{definition}[Confluence]
A topology operator is confluent if the sequence of statistical inputs
it emits depends only on the order of observations within each input,
and not on how the inputs interleave on arrival.
\end{definition}

Positional pairing is confluent, as is greedy time alignment of
time-ordered inputs. This is the determinacy property of Kahn
networks~\cite{kahn1974}. The independent topology is not confluent: its
output order is the arrival order. A key-based join that emits a row
when its last component arrives is not confluent either. For a
confluent topology the trajectory is independent of scheduling; for a
non-confluent one the interleaving is part of the experiment and the
ordering policy has to be recorded (\cref{sec:ordering}).

\section{Execution Runtime}
\label{sec:runtime}

The runtime determines how observations reach the run. An observation
may be \emph{pushed} by a source or \emph{pulled} by the runtime, and
execution may be synchronous or asynchronous. These are properties of
execution and leave \eqref{eq:transition} unchanged. In particular, a
procedure does not need an asynchronous interface because its
observations arrive asynchronously (\cref{pr:execution}).

\subsection{Concurrency and ordering}
\label{sec:ordering}

Several producers may deliver observations concurrently. The runtime may
acquire concurrently, but the inferential state of one procedure
instance has one well-defined transition order
(\cref{pr:single-writer}):
\[
\text{concurrent acquisition}
\;\longrightarrow\;
\text{ordered delivery}
\;\longrightarrow\;
\text{serialized state transition}.
\]
This is preferable to requiring every procedure to be thread-safe. It
keeps locking out of statistical code, and it makes the order in which
observations were applied a recorded fact, where otherwise it would be
an accident of scheduling.

The ordering policy names the order that is intended: arrival order,
per-source sequence-number order, event-time order, or an explicit
ordering key. A per-source policy does not fix the interleaving
\emph{between} sources. If the topology is confluent that does not
matter. If it is not, either a total order is imposed, for instance by
merging on event time, or the realized order of delivery is logged so
that it can be replayed (\cref{sec:replay}).

\subsection{Batching}
\label{sec:batching}

A runtime may deliver a batch $B=(x_{t+1},\ldots,x_{t+k})$ in one call.
Write $F^{*}$ for the iterated transition, that is, the fold of $F$ over
a sequence of inputs that stops at a terminal state as in
\cref{def:trajectory}. Batching is \emph{transparent} if
\begin{equation}
F^{*}\bigl(S_t,R_t;\,B\bigr)
\;=\;
F^{*}\bigl(S_t,R_t;\,x_{t+1},\ldots,x_{t+k}\bigr),
\label{eq:batch}
\end{equation}
including all intermediate outputs. A transparent batch is a transport
optimization and nothing else.

The condition is stricter than it looks for procedures with a stopping
rule. If a decision is evaluated only at the end of each batch, the
procedure can stop only at batch boundaries. That is a group-sequential
design~\cite{jennison2000}, a different statistical procedure with
different error probabilities and sample-size distribution. In that case
the batch structure belongs to the definition of $F$ and to the input
contract, not to the runtime.

\subsection{Duplicates and delivery guarantees}
\label{sec:delivery}

Distributed transports duplicate and reorder messages. Using source
identifiers and sequence numbers, the runtime can deduplicate, hold
observations back until a sequence is contiguous, or reject gaps
outright. Combined with a checkpoint that records the per-source
position atomically with the inferential state
(\cref{sec:persistence}), at-least-once delivery yields state
transitions that are applied exactly once: after a restart the source is
rewound to its recorded position or earlier, and anything already
applied is dropped as a duplicate.

The procedure is never responsible for broker offsets or
acknowledgement protocols.

\subsection{Backpressure}

Inputs may arrive at very different rates. The runtime can respond by
buffering, blocking producers, throttling, sampling, or dropping. Only
blocking and lossless buffering are pure execution policies. Sampling
and dropping change the observation history, so the policy and the
number of affected observations become part of the experiment and must
appear in its provenance.

\section{Runs: Lifecycle, Failure, and Persistence}
\label{sec:run}

A run is one concrete execution of a procedure. A runtime is a
disposable environment that drives it. An \emph{experiment} is a
higher-level scientific object that may contain several runs.

\subsection{Initialization and restoration}

A new run may take its initial state $S_0$ from the procedure's default,
an explicit prior, historical data, or the final state of an earlier
run. All of these \emph{initialize} a new run: it has a new identity and
an empty history. \emph{Restoring} a checkpoint is different. It
continues an existing run with its identity and history intact. Two runs
may hold numerically identical states and still differ in provenance,
and the architecture keeps the two operations apart for that reason.

\subsection{Lifecycle}
\label{sec:lifecycle}

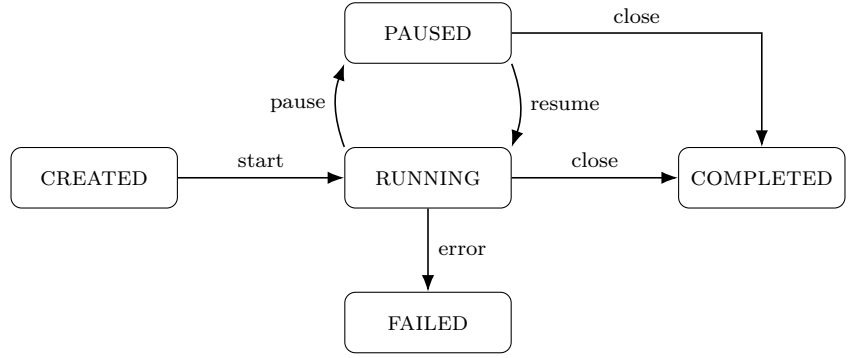
\begin{figure}[t]
\centering
\begin{tikzpicture}[node distance=11mm and 22mm]

\node[lifecycle] (created) {created};
\node[lifecycle, right=of created] (running) {running};
\node[lifecycle, right=of running] (completed) {completed};
\node[lifecycle, above=of running] (paused) {paused};
\node[lifecycle, below=of running] (failed) {failed};

\draw[arrow] (created) -- node[above,lbl] {start} (running);
\draw[arrow] (running) -- node[above,lbl] {close} (completed);
\draw[arrow] (running.north west) to[bend left=20]
    node[left,lbl] {pause} (paused.south west);
\draw[arrow] (paused.south east) to[bend left=20]
    node[right,lbl] {resume} (running.north east);
\draw[arrow] (paused.east) -| node[pos=.25,above,lbl] {close} (completed.north);
\draw[arrow] (running) -- node[right,lbl] {error} (failed);

\end{tikzpicture}
\caption{
Run lifecycle. A run restored from a checkpoint starts in
\textsc{paused}, and attaching a runtime resumes it. A failed run is
never continued in place; recovery restores an earlier checkpoint.
}
\label{fig:lifecycle}
\end{figure}

A procedure has a minimal lifecycle: it is active until its state is
terminal. A run has a richer one (\cref{fig:lifecycle}). The run owns
the execution lifecycle; the procedure owns the inferential semantics.

Reaching a terminal state of the procedure does not by itself end the
run. The runtime may notify consumers, drain or discard pending data,
close sources, and only then close the run. Inputs that arrive after a
terminal state are not applied (\cref{sec:terminal}); they are counted.

\subsection{Errors and transactional updates}
\label{sec:errors}

Failures occur at different boundaries and are the responsibility of
different layers (\cref{tab:errors}).

\begin{table}[t]
\centering
\caption{Error classes, the layer responsible for them, and their
consequence. None of them modifies the inferential state.}
\label{tab:errors}
\begin{tabular}{lll}
\toprule
\textbf{Error} & \textbf{Responsible layer} & \textbf{Consequence} \\
\midrule
Source error & adapter, runtime & run stays valid \\
Delivery or ordering error & runtime, topology & run fails \\
Validation error & contract, topology & skipped and reported, or fails \\
Numerical or inference error & procedure & transition not committed; run fails \\
Serialization error & persistence & previous checkpoint stays valid \\
Consumer error & consumer & recorded; inference continues \\
\bottomrule
\end{tabular}
\end{table}

That no failure modifies the inferential state is by design. A state
update is a transaction,
\[
\text{validate}
\;\longrightarrow\;
\text{compute}
\;\longrightarrow\;
\text{commit},
\]
and anything that fails before the commit leaves $(S_{t-1},R_{t-1})$ in
place. The simplest way to obtain this is to make states immutable
values and $F$ a function that returns a new state: there is then
nothing to roll back.

A consumer failure is not an inference failure
(\cref{pr:consumers}). A visualization that crashes must not invalidate
an otherwise valid trajectory.

An invalid input is not an inference failure either. A run does not
control the stream it consumes, and a single malformed message must not
make it unusable. The reference implementation therefore skips invalid
inputs by default: an observation the topology cannot use, a statistical
input that violates its contract, or one that the transition itself
finds invalid. Nothing is committed, so skipping is exact. Every such
input is reported as an incident, whatever the policy, and counted
when it is skipped.
Stopping the run on invalid input is an explicit choice, and the
experimenter may name the error it raises. A separate threshold stops
the run after a given number of consecutive invalid inputs: one bad
message is noise, a source that sends nothing valid is broken. The policy is recorded in
provenance. Contracts themselves are strict: a value of a foreign
numeric type, such as a \code{numpy} scalar, is refused with a message
that says how to convert it in the adapter, rather than interpreted.

\subsection{Checkpoints}
\label{sec:persistence}

We distinguish a \emph{snapshot}, which represents state for inspection,
from a \emph{checkpoint}, which is sufficient to continue the run.
\Cref{tab:checkpoint} lists what a run checkpoint contains. It contains
no sockets, threads, event loops, or other runtime resources
(\cref{pr:persistent}).

\begin{table}[t]
\centering
\caption{Contents of a run checkpoint.}
\label{tab:checkpoint}
\begin{tabular}{lp{0.62\linewidth}}
\toprule
\textbf{Component} & \textbf{Purpose} \\
\midrule
Procedure identity & name, version, and configuration; checked on restore \\
Inferential state $S_t$ & continues the mathematics \\
Random state $R_t$ & continues a randomized procedure exactly \\
Topology state & observations consumed from sources but not yet applied \\
Source positions & where each source resumes: file offsets, cursors,
broker offsets, sequence numbers \\
History digest & identifies the input history that led to $S_t$ \\
Observation log position & record count, byte offset and hash chain of
the observation log, if one is kept \\
Provenance & see \cref{sec:provenance} \\
Format version & allows the checkpoint format to evolve \\
\bottomrule
\end{tabular}
\end{table}

\paragraph{Consistency.}
A checkpoint must describe a single cut through the execution: every
observation before a source's recorded position is either reflected in
$S_t$ or held in the topology state, and no observation at or after that
position is. With a single writer this is easy to guarantee. Checkpoints
are taken between transitions, when all components refer to the same
prefix of every source. A system that runs topology operators and the
procedure in separate processes needs a distributed snapshot protocol to
obtain the same cut~\cite{chandy1985,elnozahy2002}. Restoration from a
consistent checkpoint continues the trajectory from $S_t$ without
altering the semantics of the history already processed.

\paragraph{History digest.}
Let $h_0$ be a fixed constant and
$h_t=H\bigl(h_{t-1}\,\|\,\mathrm{enc}(X_t)\bigr)$ for a cryptographic
hash $H$ and a canonical encoding of statistical inputs. Two runs with
equal digests applied the same input history, up to hash collisions. The digest lets a replay
be verified without storing the history in the checkpoint.
For the digest to be comparable across languages, the encoding must not
depend on how a language prints floating-point numbers. The reference
implementation uses the JSON Canonicalization Scheme~\cite{rfc8785},
which fixes number formatting (the shortest representation that
round-trips), string escaping and key order. Integers that a binary64
value cannot hold exactly are refused by the input contracts, and 64-bit
values such as $R_t$ are stored as strings. The same requirement applies
earlier, to observations: one whose value cannot be recorded exactly in
JSON (a non-finite number, an integer beyond $2^{53}$, a string that
UTF-8 cannot encode, a source name that is not a string) is invalid
input, rejected before it reaches a topology buffer, so that neither the observation log nor a checkpoint ever holds
a value that a replay would read back differently.

\subsection{Provenance}
\label{sec:provenance}

A reproducible run records the procedure identity, version, and
configuration; the input contracts; the identity of every source; the
topology; the ordering and delivery policy; the runtime mode of every
execution segment; the numerical backend; the seed; units and
calibration where relevant; counts of everything that was skipped or
dropped; and the lifecycle events with their timestamps. Provenance is
part of the scientific representation of a run. It is not an operational
log that can be discarded once the run has finished.

\subsection{Replay and simulation}
\label{sec:replay}

Because a procedure sees only logical inputs, the same input can be
supplied by a live source, by historical data, by a recorded
observation log, or by a simulator. Replay is then ordinary execution with a
different adapter. It supports debugging, regression testing, comparison
of methods on identical histories, and regeneration of figures.

An observation log is part of the run's persistent state. It records
the observations delivery passes on, not every arrival: duplicates and
observations still held back are left out, so that its position and the
source positions describe the same cut. The checkpoint
records its position, and restoring the checkpoint verifies the log
against the recorded hash chain and cuts it back to that position:
observations logged after the checkpoint belong to an execution that
was abandoned. The run does not preserve them. After a restore they are
delivered again by the sources from their recorded positions, and a
source that cannot deliver past observations again must be buffered by
the application, between the source and the run.

A simulator can also close the loop of an adaptive experiment,
\[
\text{procedure output}
\;\longrightarrow\;
\text{experimental design}
\;\longrightarrow\;
\text{simulator}
\;\longrightarrow\;
\text{new observation},
\]
without any change to the procedure.

\section{Numerical Semantics}
\label{sec:numerical}

Long-running procedures accumulate many small terms. Whether a run can
be reproduced, and whether two implementations agree, are therefore
questions about floating-point
arithmetic~\cite{goldberg1991,higham2002} as much as about statistics.
The usual mechanisms apply: likelihoods are kept in the log domain,
normalizing constants are computed with a stable log-sum-exp, long sums
use compensated accumulation~\cite{kahan1965,neumaier1974}, and
non-finite values are rejected before they are committed. The runtime
enforces the last point for every procedure: a transition that
overflows, or yields a non-finite state or output, raises a numerical
error and commits nothing. Procedures that estimate a variance also
report a numerical error when rescaled data make the estimate
underflow, rather than deciding on a meaningless statistic.

More important for the architecture is to say what ``equal
trajectories'' means. We distinguish four levels.

\begin{description}[leftmargin=2em,style=nextline,itemsep=2pt]
    \item[Bitwise equality.]
    The same implementation on the same platform performs the same
    operations in the same order. Checkpoint, replay, batch, and
    synchronous/asynchronous equivalence should hold at this level,
    because none of them is supposed to change which operations are
    performed.
    \item[Tolerance equality.]
    Different backends, compilers, or languages agree on states and
    continuous outputs within declared absolute and relative tolerances.
    \item[Decision equality.]
    Discrete outputs (decisions, alarms, counts) are equal. Tolerance
    equality does not imply decision equality: a decision is a
    discontinuous function of the state, and two implementations that
    agree to $10^{-12}$ can disagree about a threshold crossing. A
    numerical contract must say how such cases are treated, for example
    by reporting the margin to the boundary and flagging decisions taken
    within tolerance of it.
    \item[Distributional equality.]
    Randomized procedures driven by different generators agree only in
    distribution, unless the generator itself is part of the contract.
\end{description}

Compensation terms are part of the inferential state. A compensated sum
that is restored without its compensation term continues with a
different value than the uninterrupted run.

Bitwise equality is fragile in ways that have nothing to do with the
statistics. Python's built-in \code{sum} adds floating-point numbers with
compensation since version~3.12 and without it before, so a particle
filter whose weights are summed with it follows a different trajectory,
in the last bits, on different versions of the same language. The
reference implementation found this through its conformance vectors and
now sums floats on the inference path only by an explicit left-to-right
loop, or by the compensated sum whose correction term is stored in the
state.

The levels can be tested across implementations only if there is
something to compare against. The reference implementation ships
language-independent conformance vectors: for each of sixteen cases,
the inputs, every output, the final state, the random state and the
history digest, all encoded canonically. A second implementation
conforms bitwise if it reproduces them exactly, and at the tolerance or
decision level if it reproduces them within its declared contract.

\section{Consumers, Adapters, and Feedback}
\label{sec:consumers}

\subsection{Adapters}
\label{sec:adapters}

Physical and external systems are connected through adapters: imaging
archives speaking DICOM~\cite{dicom}, industrial equipment exposed via
OPC~UA~\cite{opcua}, sensor networks using MQTT~\cite{mqtt}, databases,
files, laboratory instruments, and network services. A medical imaging
pipeline, for instance, is
\[
\text{scanner}
\rightarrow
\text{DICOM}
\rightarrow
\text{adapter}
\rightarrow
\text{logical image input}
\rightarrow
\text{topology}
\rightarrow
\text{run}.
\]
DICOM metadata may be used for routing, pairing, filtering, and
provenance without becoming part of the statistical procedure. The same
logical input can be rebound to another scanner, an archive, or a
phantom simulator, provided the new source satisfies the input contract.

\subsection{Consumers}

A consumer receives outputs $O_t$. It has no access to the inferential
state and cannot change it (\cref{pr:consumers}). Typical consumers are
visualization, logging, reporting, database storage, alerting,
experimental design, and control or actuator integration. The first
four are passive. The last three act on the world.

\subsection{Feedback}
\label{sec:feedback}

\begin{figure}[t]
\centering
\begin{tikzpicture}[node distance=11mm and 11mm]

\node[arch] (obs) {Observations\\$X_t$};
\node[arch, right=of obs] (proc) {Sequential\\Procedure};
\node[arch, right=of proc] (out) {Output\\$O_t$};
\node[arch, right=of out] (cons) {Consumer};
\node[arch, below=of cons] (act) {Action $A_t$};
\node[arch, below=of obs] (world) {External\\system};

\draw[arrow] (obs) -- (proc);
\draw[arrow] (proc) -- (out);
\draw[arrow] (out) -- (cons);
\draw[arrow] (cons) -- (act);
\draw[arrow] (act) -- node[above,lbl] {actuation} (world);
\draw[arrow] (world) -- node[right,lbl,align=left]
    {response and\\applied action} (obs);

\end{tikzpicture}
\caption{
Inference and action are separated. The loop closes through the external
system. An action that influences later observations is itself recorded.
}
\label{fig:feedback}
\end{figure}
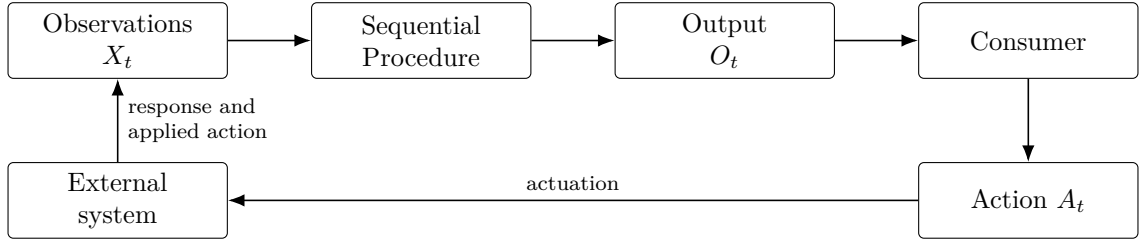

Inference is $X_{1:t}\mapsto S_t\mapsto O_t$. A consumer may implement
$O_t\mapsto A_t$, where $A_t$ is an external action. If the action
changes the subsequent observation process, so that $X_{t+1}$ depends on
$A_t$, the system contains a feedback loop (\cref{fig:feedback}).

The separation lets one procedure serve passive and active experiments
alike. It leaves two obligations, which should be stated plainly.

\begin{remark}[Separation does not confer validity]
Moving actuation out of the procedure does not make inference valid
under feedback. Whether a procedure remains valid when the design
depends on earlier outputs is a statistical property. Likelihood-ratio
and e-process arguments survive designs that are predictable with
respect to the past~\cite{ramdas2023}; many fixed-sample arguments do
not. The architecture contributes two things here. It makes the
dependence explicit, and it lets the same $F$ be analyzed with and
without the loop.
\end{remark}

\begin{remark}[Actions are data]
An action that changes the observation process is part of the history
of the experiment. It should either re-enter the statistical boundary as
a declared logical input, as a control input does in a Kalman filter, or
be recorded in provenance. A consumer never calls back into the run. The
loop always closes through a source.
\end{remark}

\section{Composition and Stateful Predictive Models}
\label{sec:pipeline}

\subsection{Composition}

Procedures compose with one another and with stateless transformations.
In
\[
\text{raw input}
\rightarrow
\text{exponential moving average}
\rightarrow
\text{residual}
\rightarrow
\text{change-point procedure}
\]
the moving average is a stateful transformation and therefore itself a
procedure, the residual is a stateless map, and the detector is a
procedure. Feeding an e-process into a CUSUM-type detector is another
composition of the same shape.

Composition is closed: the chain of two procedures is a procedure whose
state is the pair of their states and whose transition applies them in
order. No component inherits from a pipeline type. Composition is an
external mechanism and adds no depth to any hierarchy
(\cref{pr:portability}).

\subsection{Stateful predictive models}
\label{sec:ml}

The abstraction is not restricted to classical statistics. A recurrent
model with state $S_t=(h_t,c_t)$ and update $(S_t,O_t)=F(S_{t-1},X_t)$,
such as an LSTM~\cite{hochreiter1997} used for prediction, satisfies
\cref{def:procedure}. Its trained weights are configuration; they are
part of the procedure's identity and do not change during a run. If the
weights are updated online, they become part of the state.

A stateless mapping $O=f(X)$ is not a procedure. It is a transformation
in the topology or pipeline. This distinction admits stateful predictive
models without turning the architecture into a general machine-learning
framework. The criterion for membership is whether a method maintains
state and produces outputs sequentially. Which statistical school it
belongs to is irrelevant.

\section{Reference Procedure: The Kiefer--Weiss Problem}
\label{sec:kw}

We use a sequential test as the reference procedure because it
exercises most of the architecture at once: log-domain accumulation, a
stopping boundary, a terminal decision, and state that must survive
long runs.

\paragraph{Problem.}
Let $X_1,X_2,\ldots$ be independent with density $f_\theta$, and
consider testing $\theta_0$ against $\theta_1>\theta_0$ with error
probabilities at most $\alpha_0$ and $\alpha_1$. The SPRT minimizes the
expected sample size at $\theta_0$ and at
$\theta_1$~\cite{waldwolfowitz1948}. Between them its expected sample
size can be much larger, and its sample size is unbounded. Kiefer and
Weiss~\cite{kiefer1957} asked for the test that minimizes the
\emph{maximum} expected sample size $\sup_\theta \Exp_\theta N$ under the
same error constraints. The \emph{modified} Kiefer--Weiss problem
minimizes $\Exp_{\theta^*}N$ at one intermediate point $\theta^*$; in
symmetric normal and binomial problems the two
coincide~\cite{weiss1962}.

\paragraph{The 2-SPRT.}
Lorden~\cite{lorden1976} proposed running two one-sided SPRTs against
the intermediate point. With
\begin{equation}
L^{0}_t=\sum_{i=1}^{t}\log\frac{f_{\theta^*}(X_i)}{f_{\theta_0}(X_i)},
\qquad
L^{1}_t=\sum_{i=1}^{t}\log\frac{f_{\theta^*}(X_i)}{f_{\theta_1}(X_i)},
\label{eq:two-sprt}
\end{equation}
the test stops at $N=\inf\{t: L^0_t\ge a_0 \text{ or } L^1_t\ge a_1\}$
and rejects $\theta_0$ if the first boundary was reached and $\theta_1$
if the second was. Lorden showed that the 2-SPRT approximately minimizes
$\Exp_{\theta^*}N$ among tests with the same error probabilities, and
Huffman~\cite{huffman1983} gave an explicit choice of $\theta^*$ for
which it is an asymptotically efficient solution of the original minimax
problem.

Under $\theta_i$ the process $\exp(L^i_t)$ is a non-negative martingale
with unit mean, so by Ville's inequality~\cite{ville1939}
\begin{equation}
\Prob_{\theta_0}(\text{reject }\theta_0)\le e^{-a_0},
\qquad
\Prob_{\theta_1}(\text{reject }\theta_1)\le e^{-a_1}.
\label{eq:ville}
\end{equation}
The thresholds $a_i=\log(1/\alpha_i)$ therefore guarantee the error
rates. They are conservative, because they ignore that the other
boundary may be reached first.

\paragraph{Choice of $\theta^*$.}
To first order the two boundaries should be reached at the same time
under $\theta^*$. With $I(\theta,\theta')$ the Kullback--Leibler
divergence, this gives
\begin{equation}
\frac{a_0}{I(\theta^*,\theta_0)}=\frac{a_1}{I(\theta^*,\theta_1)} .
\label{eq:theta-star}
\end{equation}
For normal observations with known variance $\sigma^2$ the solution is
$\theta^*=\theta_0+(\theta_1-\theta_0)\sqrt{a_0}/(\sqrt{a_0}+\sqrt{a_1})$,
the midpoint when $a_0=a_1$. In this case the continuation region is a
triangle in the $(t,\sum_i X_i)$ plane, and the test always stops by
\begin{equation}
n_{\max}=\left\lceil
\frac{2\sigma^2\,(a_0/\delta_0+a_1/\delta_1)}{\delta_0+\delta_1}
\right\rceil,
\qquad \delta_i=|\theta^*-\theta_i| .
\label{eq:nmax}
\end{equation}

\paragraph{As a Sequential Procedure.}
The state is $S_t=(t,L^0_t,L^1_t,d_t)$, with each sum carried together
with its compensation term and $d_t$ the decision, if any. The terminal
set is $\{d_t\neq\varnothing\}$. The output repeats the statistics and
the decision. The procedure owns the meaning of this state and its
stopping semantics. It owns no broker integration, scheduling,
acquisition, checkpoint storage, visualization, or reporting.

\paragraph{Executing optimal plans.}
The 2-SPRT is an approximation. Exact solutions of the Kiefer--Weiss
problem can be computed numerically by backward induction for Bernoulli
observations, discrete exponential families, and normal
observations~\cite{novikov2022a,novikov2022b,novikov2023}. The result of
such a computation is a sampling plan: a horizon $H$ and, for every step
$n\le H$, an interval of the running sum $\sum_{i\le n}X_i$ in which
sampling continues. Design and execution are different problems. The
design is computed once, offline. Its execution on a stream is a
Sequential Procedure with state $(n,\sum_{i\le n}X_i,d_n)$, and the
reference implementation provides one that takes any plan of this form
as its configuration. The plan thereby becomes part of the procedure's
identity and is recorded in provenance. As a check, the plan formed by
the two boundary lines of the Gaussian 2-SPRT reproduced the stopping
time and decision of the 2-SPRT on all $4\,000$ simulated paths we
compared.

The reference implementation also contains the design side, following
\cite{novikov2022a,novikov2022b,novikov2023}: backward induction
on the Lagrangian $E_{\theta^*}N+\lambda_0P_{\theta_0}(\text{reject})+
\lambda_1P_{\theta_1}(\text{accept})$ over plans with horizon $H$
\cite{lorden1980}, with the multipliers adjusted until the error
probabilities meet their targets. The running sum is carried on a
lattice: the integers for Bernoulli and Poisson data, where every
quantity is exact, a lattice of step $h$ for normal data, and the
midpoints of cells of width $h$ for exponential data, where the error is
of order $h^2$. Writing the error probabilities as expectations
under $\theta^*$ of the likelihood ratios lets one backward pass return
the Lagrangian and both error probabilities at once, so that adjusting
the multipliers by Broyden's method takes a few passes. The resulting
plan is executed by the same plan executor. For $\theta_0=0.3$ against
$\theta_1=0.5$ with $\alpha_0=\alpha_1=0.05$ and $H=150$, it attains
error probabilities $0.0498$ and $0.0496$ and a maximum expected sample
size of $45.2$, against $50.8$ for Wald's SPRT with thresholds
calibrated to the same error probabilities ($56.9$ with Wald's
approximate thresholds $\pm\log 19$). The value of the Lagrangian obtained by backward
induction agrees with the forward evaluation of the plan to machine
precision, and the exact operating characteristic agrees with
simulation through the executor. For a normal mean, $\theta_0=0$ against
$\theta_1=0.5$ with $\sigma=1$ and $H=120$, the plan computed on a lattice of step
$\sigma/16$ has error probabilities $0.0515\pm0.0015$ when executed on
simulated normal data, and an expected sample size at $\theta^*$ of
$31.0$. The 2-SPRT with thresholds calibrated to the same error
probabilities ($a=2.16$) needs $31.3$: in this example Lorden's
approximation is within about one percent of the optimum, which the
backward induction makes it possible to quantify.

The modified problem fixes $\theta^*$. The Kiefer--Weiss problem asks
for the $\theta^*$ at which the resulting test has the smallest maximum
expected sample size, and by Lorden's characterisation~\cite{lorden1980}
that is the $\theta^*$ at which the maximum is attained at $\theta^*$
itself. For discrete data the implementation generates candidates by
bisection on $g(\theta^*)=\arg\max_\theta E_\theta N-\theta^*$, designing
each in full, and returns the candidate, the first-order point included,
with the smallest maximum expected sample size. For the Bernoulli example the first-order
$\theta^*=0.3971$ gives a maximum expected sample size of $45.17$,
attained at $0.3933$; the search returns $\theta^*=0.3925$ with $45.01$,
attained at $0.3936$. Near the solution the maximum is flat in
$\theta^*$, and for discrete data it moves in small steps because the
attainable error probabilities do, so the gain over the first-order
point is small and the characterisation holds only approximately.

For normal data the characterisation itself is well determined: the
function $g(\theta^*)=\arg\max_\theta E_\theta N-\theta^*$ is smooth and
decreasing, and on a symmetric problem its root is the midpoint to
$3\cdot10^{-4}$ of $\theta_1-\theta_0$. The maximum expected sample size
is not: over a tenth of the interval around the root it changes by about
$2\cdot10^{-4}$ of its value, the order of the calibration error of the
multipliers, so that which candidate has the smallest maximum would be
partly decided by rounding, and would change with the candidates a search
happens to visit. For $\theta_0=0$ against $\theta_1=0.5$ with error
probabilities $0.05$ and $0.2$ and $H=60$, the candidate with the smallest
maximum among those of a bisection over the whole interval is
$\theta^*=0.340$ with $17.131$, the root of the characterisation is
$0.3462$ with $17.140$, and the first-order point $0.289$ gives $17.148$:
the three plans differ by less than a tenth of a percent. For data on a
lattice (normal, exponential) the implementation therefore returns the
root of the characterisation, found by the Illinois variant of regula
falsi to $2\cdot10^{-4}$ of $\theta_1-\theta_0$; it is reproducible to that
tolerance, and a symmetric problem gets the midpoint to it.

Because the root is well determined, it can be tabulated. In the
standardised coordinate $\lambda=(\theta^*-\theta_0)/(\theta_1-\theta_0)$
it depends on $\alpha_0$, $\alpha_1$ and the horizon, which enters
through $\rho=\delta\sqrt H/(z_{\alpha_0}+z_{\alpha_1})$, the ratio of the
horizon to the shortest fixed-sample test; exchanging the hypotheses maps
$\lambda$ to $1-\lambda$, so the correction to the first-order point is
antisymmetric in $(\alpha_0,\alpha_1)$. A pilot study showed it smooth in
$\log\alpha_0$, $\log\alpha_1$, nearly constant in $\rho$ above $1.5$
(and in $\delta$ at fixed $\rho$, to $1.5\cdot10^{-3}$), with a boundary
layer of up to $0.015$ as $\rho$ approaches $1$. The roots were computed
on Gauss--Lobatto--Legendre nodes, which are D-optimal for polynomial
regression~\cite{guest1958}, in tensor product (63 nodes for $\alpha_0,\alpha_1\in[0.005,
0.2]$, $\rho\in[1.5,2.5]$, only on one side of the diagonal), and fitted
by an antisymmetric Legendre series of degree $5$ in each $\log\alpha$
and $1$ in $\rho$; on 11 Halton points not used in the fit the error is
$1.2\cdot10^{-4}$ rms and at most $3.3\cdot10^{-4}$. The search for the
least favourable $\theta^*$ starts from a bracket of $\pm0.004$ around
the tabulated root ($\pm0.03$ below $\rho=1.5$) and falls back to the
whole interval when the bracket does not contain a sign change. The
example above takes 83 seconds; the bisection over the whole interval
with the earlier implementation, without numpy and without the
calibration below, took two hours.

Each candidate is designed in full. The backward pass is written as
array operations (numpy, when installed), and once Broyden's method has
come within $10^{-3}$ of the targets the multipliers are adjusted by a
small joint raise before the slower coordinate search is tried; on
thirteen designs across the four families both changes leave the plans
unchanged.

For discrete data there is a further limit. The plans that are optimal
for some pair of multipliers are non-randomised, and near the smallest
feasible horizon they need not reach every pair of error probabilities
that some test of that size reaches; randomised boundaries would. The
implementation then tries smaller horizons, which also respect the
constraint, and reports a failure only if none succeeds. Before any of
this it compares the targets with the Neyman--Pearson bound: a plan
that stops by $H$ is a test based on $H$ observations, so its
$\alpha_1$ is at least that of the most powerful randomised level-$\alpha_0$
test of that size, computed exactly for Bernoulli, Poisson and normal
data. Targets below the bound are reported at once, instead of after a
search that cannot succeed.

\paragraph{Validation.}
A reference procedure can be validated against analytical facts without
any other implementation. For the 2-SPRT these are the error bound
\eqref{eq:ville}, the closed-form $\theta^*$, and the sample-size bound
\eqref{eq:nmax}, which constant data equal to $\theta^*$ attains exactly.
\Cref{tab:oc} adds a Monte Carlo study for $\theta_0=0$, $\theta_1=0.5$,
$\sigma=1$, and nominal error rates $0.05$. The SPRT with Wald's
thresholds attains error rates of about $0.038$. The 2-SPRT with the
guaranteed thresholds is more conservative and, as a result, slower. Once
its threshold is calibrated to the SPRT's attained error rate, the
trade-off the Kiefer--Weiss problem is about becomes visible. The
expected sample size at the midpoint drops from $42.4$ to $37.1$, the
expected sample sizes at the hypotheses rise from about $24.3$ to about
$27.0$, and the largest sample observed falls from $384$ to $77$,
within the bound $n_{\max}=79$.

\begin{table}[htb]
\centering
\caption{
Operating characteristics for normal observations with unit variance,
$\theta_0=0$ against $\theta_1=0.5$; $20\,000$ replications per cell,
standard errors in parentheses. All three tests see the same simulated
data at each $\theta$. The matched threshold was calibrated on a
separate random stream.
}
\label{tab:oc}
\begin{tabular}{llrrr}
\toprule
\textbf{Test} & $\theta$ & $\Exp_\theta N$ & $\max N$ & $\Prob_\theta(\text{reject }\theta_0)$ \\
\midrule
SPRT, Wald thresholds $\pm 2.944$
    & $0$    & $24.21\ (0.12)$ & $177$ & $0.0379$ \\
    & $0.25$ & $42.36\ (0.24)$ & $384$ & $0.5029$ \\
    & $0.5$  & $24.44\ (0.12)$ & $182$ & $0.9616$ \\
\addlinespace
2-SPRT, $a=\log 20\approx 2.996$
    & $0$    & $32.96\ (0.10)$ & $91$ & $0.0218$ \\
\quad ($n_{\max}=96$)
    & $0.25$ & $48.15\ (0.12)$ & $93$ & $0.5041$ \\
    & $0.5$  & $33.15\ (0.10)$ & $93$ & $0.9788$ \\
\addlinespace
2-SPRT, matched $a=2.456$
    & $0$    & $26.90\ (0.08)$ & $75$ & $0.0355$ \\
\quad ($n_{\max}=79$)
    & $0.25$ & $37.09\ (0.10)$ & $77$ & $0.5030$ \\
    & $0.5$  & $27.05\ (0.08)$ & $72$ & $0.9634$ \\
\bottomrule
\end{tabular}
\end{table}

The study drives the procedures with a bare loop over simulated values.
It uses no run, runtime, or topology. That the same transition function
serves a Monte Carlo study and a laboratory deployment unchanged is the
practical content of \cref{pr:execution}.

\section{Conformance Invariants}
\label{sec:invariants}

An architecture of this kind is not evaluated well by throughput
benchmarks. Its claims are about semantics, and they can be stated as
invariants that an implementation either satisfies or does not. Each
invariant below compares trajectories in the sense of
\cref{def:trajectory}, at the bitwise level of \cref{sec:numerical}
unless noted otherwise.

\begin{invariant}[Synchronous/asynchronous equivalence]
\label{inv:async}
For the same ordered statistical input history, synchronous and
asynchronous execution yield the same trajectory. For a confluent
topology the same holds for any interleaving of the sources.
\end{invariant}

\begin{invariant}[Checkpoint equivalence]
\label{inv:checkpoint}
For every $k$, processing $x_{1:k}$, taking a checkpoint, restoring it
from storage with no other state carried over, and processing
$x_{k+1:T}$ yields the same trajectory as processing $x_{1:T}$ without
interruption. The same holds if the
sources redeliver observations from before the checkpoint.
\end{invariant}

\begin{invariant}[Replay equivalence]
\label{inv:replay}
Replaying the recorded delivered history of a run through the same
delivery policy and topology reproduces its trajectory and its history
digest.
\end{invariant}

\begin{invariant}[Batch transparency]
\label{inv:batch}
Delivering a history in batches of any size yields the same trajectory
as delivering it element by element, including the step at which a
terminal state is reached.
\end{invariant}

\begin{invariant}[Consumer non-interference]
\label{inv:consumers}
Attaching, detaching, or crashing consumers does not change the
trajectory.
\end{invariant}

\begin{invariant}[Runtime independence]
\label{inv:runtime}
Different acquisition mechanisms that produce the same ordered logical
input history yield the same trajectory.
\end{invariant}

\begin{invariant}[Topology semantics]
\label{inv:topology}
Each topology operator reproduces known reference histories, and
operators declared confluent emit the same sequence for every
interleaving of their inputs.
\end{invariant}

\Cref{tab:invariants} relates the invariants to the principles. The
invariants are stated for a single implementation. Agreement
\emph{between} implementations in different languages is a separate
validation criterion, at the tolerance and decision levels.

\begin{table}[t]
\centering
\caption{Invariants and the principles they test.}
\label{tab:invariants}
\begin{tabular}{lll}
\toprule
\textbf{Invariant} & \textbf{Principles} & \textbf{What a violation indicates} \\
\midrule
\ref{inv:async} Sync/async & \ref{pr:execution}, \ref{pr:single-writer} & scheduling leaks into inference \\
\ref{inv:checkpoint} Checkpoint & \ref{pr:persistent} & state hidden outside the checkpoint \\
\ref{inv:replay} Replay & \ref{pr:logical}, \ref{pr:persistent} & unrecorded inputs or randomness \\
\ref{inv:batch} Batch & \ref{pr:execution} & transport changes the statistics \\
\ref{inv:consumers} Consumers & \ref{pr:consumers} & outputs feed back into state \\
\ref{inv:runtime} Runtime & \ref{pr:execution}, \ref{pr:logical} & the procedure depends on a source \\
\ref{inv:topology} Topology & \ref{pr:topology} & pairing logic inside the procedure \\
\bottomrule
\end{tabular}
\end{table}

\section{Reference Implementation}
\label{sec:implementation}

The architecture is defined semantically (\cref{pr:portability}). An
implementation must provide equivalent contracts for state transition,
input validation, output generation, lifecycle, persistence, ordering,
and reproducibility. It need not use classes or inheritance. A language
with structures and functions, traits, or type classes can express the
same contracts, and implementations in different languages need not
have the same feature set.

We provide a reference implementation in Python. It uses only the
standard library (numpy, when installed, accelerates the design of
optimal plans) and comprises about $6\,200$ lines, with one module per
layer of \cref{fig:architecture}. The contract a procedure has to
satisfy is small:

\begin{verbatim}
    config()                 -> parameters that determine F
    inputs()                 -> {name: input contract}
    initial_state()          -> S_0
    step(state, x, rng)      -> (new state, output)
    is_terminal(state)       -> bool
    encode_state(state), decode_state(value)
\end{verbatim}

States are immutable values and \code{step} is a pure function of its
arguments. The run draws the statistical random state into a fresh
generator before each step and commits it together with the new state,
which gives transactional updates without rollback logic. The random
state is a single 64-bit integer~\cite{steele2014}, so it is trivial to
checkpoint and to reproduce in another language. A checkpoint is a
plain JSON value.

The implementation includes fourteen procedures. For testing simple
hypotheses: Wald's SPRT, the 2-SPRT of \cref{sec:kw}, and an executor for
precomputed sampling plans. For clinical trials: a group sequential
test with efficacy and non-binding futility boundaries, for known or
estimated variance. For anytime-valid inference: the mixture SPRT for a
normal mean~\cite{robbins1970,johari2022}, its counterpart for unknown
variance, whose Bayes factor with a scale-invariant prior on $\sigma$ is
a test martingale for every $\sigma$~\cite{lai1976,perezortiz2024}, both
with an always-valid $p$-value and a confidence sequence, and two
nonparametric
tests by betting~\cite{shafer2021,waudbysmith2024,ramdas2023}, one for the
mean of bounded observations and a sign test of
$P(X>m)=P(X<m)$, both valid for every distribution satisfying the null
hypothesis. The sign test is a test of the median only when the
distribution has no atom at $m$: with ties at $m$, the median can be
$m$ while the signs are unbalanced, and the test then rejects. For change-point
detection: CUSUM~\cite{page1954} and the Shiryaev--Roberts
procedure~\cite{shiryaev1963,roberts1966,pollak1985}. For filtering and
estimation: a Kalman filter for the local-level model with an optional
control input, a bootstrap particle filter for the same model, a
two-sample mean difference that accepts unpaired inputs, and an
exponential moving average. Likelihood-ratio procedures take one of four
families (normal with known variance, Bernoulli, Poisson, exponential).
A design module computes operating characteristics and optimal plans
for Bernoulli, Poisson, normal and exponential data, for a fixed or the
least favourable $\theta^*$ (\cref{sec:kw}). Group sequential
boundaries are computed from alpha- and beta-spending
functions~\cite{lan1983,jennison2000} by recursive numerical
integration~\cite{armitage1969}; the O'Brien--Fleming-type efficacy
boundaries for three and five analyses reproduce the published values to
three decimals, and with estimated variance each boundary is converted to
a $t$ boundary of the same significance level. The $t$ distribution is
computed from the standard library alone: the incomplete beta function
by its continued fraction, in logarithms where the tail underflows, and
Fisher's expansion to third order when the degrees of freedom are
large; on a grid of degrees of freedom from $0.1$ to $10^{200}$ and
arguments up to $10^{300}$ it agrees with 80-digit reference values to
a relative error below $10^{-10}$. The
implementation also includes four topology operators (independent,
positional pairing, key join, time alignment), three delivery policies,
a synchronous pull runtime, and an asynchronous push runtime with
blocking or dropping backpressure. Transformations applied by a topology
or between the stages of a chain carry a name, a version and a
configuration, like procedures, and that identity is part of the
checkpoint. A checkpoint is validated field by field before any of it is
trusted, so that a damaged or foreign one is reported as such.

The invariants of \cref{sec:invariants} are executable tests.
Checkpoint equivalence, batch transparency, and consumer
non-interference are checked for fifteen combinations of procedure and
topology, including the randomized particle filter and topologies whose
buffers are non-empty at the checkpoint. The remaining invariants are
checked on representative cases. All seven hold bitwise. The suite has
$219$ tests in total; the remainder cover transactional updates,
lifecycle transitions, delivery policies, invalid-input policies,
compatibility checks on restore, the observation log, and the
statistical properties claimed for each procedure: error rates under
continuous monitoring, coverage of confidence sequences, run lengths to
false alarm, optimality of designed plans, Lorden's characterisation of
the least favourable $\theta^*$, the alpha and beta spent by group
sequential designs, and the power of designs with futility
boundaries. \Cref{tab:oc} is produced
by a script in the repository from a fixed seed.

An independent implementation of the history digest and of the SPRT in
JavaScript reproduces the conformance vectors bit for bit: all sixteen
digests, the canonical encodings of 500 numbers, and every SPRT output
for normal, Bernoulli and Poisson data.

Two examples show what the separation buys. A paired test fed by two
concurrent producers with random delays produces the same history
digest as the synchronous run on the same data, with no asynchronous
code in the procedure. A closed-loop example uses a controller that
is a consumer and acts on a simulated plant. The applied action returns
as a logical input paired with the measurement it influenced, and the
filter remains an ordinary Kalman filter with a control term.

\section{Discussion}
\label{sec:discussion}

\paragraph{What the separation provides.}
A procedure becomes reusable because its definition is no longer tied
to acquisition infrastructure. A topology becomes explicit because
pairing and alignment are represented independently of the statistical
update. A run becomes reproducible because inferential state, random
state, topology state, source positions, and provenance are persisted
together and nothing else is needed to continue. Visualization and
experimental design obtain a natural boundary: neither has to live
inside the inference code.

\paragraph{What it costs.}
The architecture asks more of a method author in one respect. The state
must be explicit and serializable, and all randomness must pass through
$R_t$. Methods written as loops over arrays have to be rewritten as
transitions. This is a small change for recursive statistics and a
larger one for methods that revisit past data. Such
methods must carry the data they revisit in their state, at the
corresponding cost in checkpoint size.

\paragraph{Limitations.}
The reference implementation covers fourteen procedures; optimal
design is limited to simple hypotheses about one-parameter families, and
continuous data are designed on a lattice. It demonstrates the contract and the
conformance tests; it is not yet the general-purpose library that the
gap described in \cref{sec:introduction} calls for. It is a single
process and has not been evaluated for throughput. Its adapters read
CSV files, simulate, and replay recorded observation logs; it contains
none for instruments or message brokers. Cross-language agreement is
demonstrated bitwise for the history digest and for the SPRT, but not
yet for a complete second implementation, nor at the tolerance and
decision levels for procedures that use transcendental functions, whose
last bits may differ between mathematical libraries. Finally, the architecture
does not address the statistical validity of adaptive experiments
(\cref{sec:feedback}). It makes the adaptation visible, and establishing
validity remains the job of the statistical analysis.

\paragraph{Scope.}
The architecture does not attempt to solve distributed stream
processing, machine learning, laboratory automation, or control. Those
systems integrate with it through the boundaries described above. The
restriction is intentional: the objective is a stable semantic boundary
around sequential scientific inference, small enough to be implemented
and verified.

\section{Conclusion}
\label{sec:conclusion}

We proposed an architecture for sequential scientific analysis in which
statistical inference is separated from acquisition, topology,
execution, persistence, and external action. The central abstraction is
the Sequential Procedure,
\[
(S_t,\,O_t,\,R_t) \;=\; F(S_{t-1},\,X_t,\,R_{t-1}),
\]
where $X_t$ may carry any non-empty subset of several logically distinct
inputs. Data topology determines how inputs are related. The runtime
determines how they are acquired and ordered. A run owns state,
lifecycle, persistence, and provenance. Adapters bind logical inputs to
physical systems. Consumers interpret outputs and may act on the world:

\begin{quote}
\textbf{Procedure infers; consumers act.}
\end{quote}

The separation is testable. Seven invariants follow from the principles,
and a reference implementation satisfies all of them bitwise for
procedures ranging from the 2-SPRT for the Kiefer--Weiss problem to a
randomized particle filter. The architectural claim is not tied to a
statistical method or a programming language: the state transition of a
sequential scientific procedure should remain independent of the
mechanisms that acquire observations, execute computation, persist
state, and act on outputs. That independence is what makes sequential
experiments reproducible, composable, and long-lived.

\section*{Code availability}

The reference implementation is released as the Python package
\texttt{sequential-analysis}, version~0.1.0, available from the Python
Package Index (\url{https://pypi.org/project/sequential-analysis/}) and
at \url{https://github.com/andre-a-no/sequential-analysis}, together with
the conformance vectors, the JavaScript check and the scripts that
produce \cref{tab:oc} and the table of least favourable parameters. It
is distributed under the GNU Affero General Public License, version~3,
or under a commercial license from Novikov Laboratories LLC.

\printbibliography

\end{document}